\input amstex

\documentstyle{amsppt}
\loadbold

\def\<{\left<}										
\def\>{\right>}

\def\rank{\text{rank}}

\nologo

\def\qedd{
\hfill
\vrule height4pt width3pt depth2pt
\vskip .5cm
}

\magnification=\magstephalf

\topmatter
\title
The Dirac operator of a commuting $d$-tuple
\endtitle

\author William Arveson
\endauthor

\affil Department of Mathematics\\
University of California\\Berkeley CA 94720, USA
\endaffil

\date 11 June, 2000
\enddate
\thanks On appointment as a Miller Research 
Professor in the Miller Institute for Basic 
Research in Science.  
Support is also acknowledged from 
NSF grant DMS-9802474
\endthanks
%
%
%
\abstract 
Given a commuting $d$-tuple $\bar T=(T_1,\dots,T_d)$ 
of otherwise arbitrary nonnormal 
operators on a Hilbert space, there is an 
associated Dirac operator $D_{\bar T}$.  Significant attributes 
of the $d$-tuple are best expressed in terms of $D_{\bar T}$, 
including the Taylor spectrum and the notion of 
Fredholmness.  

In fact, {\it all} 
properties of $\bar T$ derive from its Dirac operator.  We 
introduce a general notion of Dirac operator (in 
dimension $d=1,2,\dots$) that is appropriate 
for multivariable operator theory.  We show that every abstract 
Dirac operator is associated with a commuting $d$-tuple, and 
that two Dirac operators are isomorphic iff their associated 
operator $d$-tuples are unitarily equivalent.  

By relating the curvature invariant introduced in a previous paper
to the index of a Dirac 
operator, we establish a stability result for 
the curvature invariant  for
pure $d$-contractions of finite rank.  
It is shown that for the subcategory of all such $\bar T$ 
which are  a) Fredholm and and b) graded, the curvature 
invariant $K(\bar T)$ is stable under compact perturbations.  We 
do not know if this stability persists when 
$\bar T$ is Fredholm but ungraded, though there is concrete evidence that it does.     
\endabstract

\rightheadtext{Dirac operators}
\endtopmatter

\document

\subhead{Introduction}
\endsubhead
We introduce an abstract notion of Dirac operator in complex dimension 
$d=1,2,\dots$ and we show that this theory of Dirac operators 
actually coincides with the theory of commuting $d$-tuples of 
operators on a common Hilbert space $H$.  The homology 
and cohomology of Dirac operators is discussed in general terms, 
and we relate the homological picture to classical 
spectral theory by describing its application to concrete problems 
involving the solution of linear equations of the form
$$
T_1 x_1+T_2x_2+\dots+T_d x_d=y
$$
given $y$ and several commuting operators $T_1,T_2,\dots,T_d$.  

These developments grew out of an attempt to understand the 
stability properties of a curvature invariant 
introduced in a previous paper (see \cite{3}, \cite{4}), and 
to find an appropriate formula 
that expresses the curvature invariant as the 
index of some operator.  The results are presented in section
4 (see Theorem B and its corollary).  

While there is a large literature concerning Taylor's 
cohomological notion of joint spectrum for commuting 
sets of operators on a Banach space, less attention has 
been devoted to the Dirac operator 
that emerges naturally in the context of Hilbert spaces 
(however, see sections 4 through 6 of \cite{6}, where 
the operator $B + B^*$ is explicitly related to Taylor 
invertibility and the Fredholm 
property).  
We have made no attempt 
to compile a comprehensive list of references concerning the 
Taylor spectrum, but we do call the 
reader's attention to work of Albrecht \cite{1}, 
Curto \cite{5},\cite{6}, McIntosh and Pryde \cite{10}, 
Putinar \cite{12},\cite{13}, and Vasilescu \cite{15},\cite{16}.   
A more extensive list of references 
can be found in the survey \cite{6}.  Finally, I want to thank 
Stephen Parrott for useful remarks based on a draft of this paper, 
and Hendrik Lenstra for patiently enlightening me on homological issues.  

\subhead{1.  Preliminaries: Clifford structures and the CARs in dimension $d$}
\endsubhead

Since there is significant variation in the notation commonly 
used for Clifford algebras and CAR algebras, we begin with 
explicit statements of notation and terminology as it will 
be used below.  

Let $H$ be a complex Hilbert space and let $d$ be a positive integer.  
By a Clifford structure on $H$ (of 
real dimension $2d$)
we mean a real-linear mapping $R: \Bbb C^d\to\Cal B(H)$ 
of the $2d$-dimensional real vector space $\Bbb C^d$ into 
the space of self adjoint operators on $H$ which satisfies 
$$
R(z)^2 = \|z\|^2\bold 1, \qquad z\in \Bbb C^d, \tag{1.1}
$$
where for a $d$ tuple $z=(z_1,\dots,z_d)$ of complex 
numbers, $\|z\|$ denotes the Euclidean norm 
$$
\|z\|^2=|z_1|^2+\dots+ |z_d|^2.  
$$
Clifford structures can also be defined as real-linear 
maps $R^\prime$ of $\Bbb C^d$ into the space of skew-adjoint 
operators on $H$ which satisfy $R^\prime(z)^2=-\|z\|^2\bold 1$, 
and perhaps this is a more common formulation.  Note however 
that such a structure corresponds to a Clifford structure $R$ satisfying 
(1.1) by way of $R^\prime(z)=iR(z)$.  

Letting $e_1=(1,0,\dots,0),\dots,e_d=(0,\dots,0,1)$ 
be the usual unit vectors in $\Bbb C^d$
we define operators $p_1,\dots,p_d,q_1,\dots,q_d\in \Cal B(H)$ 
by $p_k=R(e_k)$, $q_k = R(ie_k)$, $k=1,\dots,d$.  The 
$2d$ operators $(r_1,\dots,r_{2d})=(p_1,\dots,p_d,q_1,\dots,q_d)$ 
are self adjoint, they  satisfy 
$$
r_k r_j + r_j r_k = 2\delta_{jk}\bold 1,
\qquad 1\leq j,k\leq 2d,  \tag{1.2}
$$
and the complex algebra 
they generate is a $C^*$-algebra isomorphic to $M_{2^d}(\Bbb C)$.  

While Clifford structures are real-linear maps of $\Bbb C^d$ there is 
an obvious way to complexify them, and once that is done one obtains 
a (complex-linear) 
representation of the canonical anticommutation relations.  This  
sets up a bijective correspondence between Clifford structures and 
reprsentations of the anticommutation relations.  The 
details are as follows.  Since 
the $2d$-dimensional real vector space $\Bbb C^d$ comes 
with an {\it a priori} complex structure, 
any real-linear mapping $R$ of $\Bbb C^d$ into the 
self adjoint operators of $\Cal B(H)$ is the real part of 
a unique complex-linear mapping $C:\Bbb C^d\to\Cal B(H)$ 
in the sense that 
$$
R(z) = C(z) + C(z)^*, \qquad z\in \Bbb C^d, \tag{1.3}
$$
and  $C$ is given by 
$C(z) = \frac{1}{2}(R(z)-iR(iz))$, $z\in \Bbb C^d$.  
Corresponding to (1.2) one finds that the operators 
$c_k= C(e_k) = \frac{1}{2}(p_k - iq_k) $, $1\leq k\leq d$ 
satisfy the canonical anticommutation 
relations 
$$
\align
c_kc_j+c_jc_k &= 0\\
c_k^*c_j+c_jc_k^*&=\delta_{jk}\bold 1. \tag{1.4}
\endalign
$$
Equivalently, the complex linear map 
$C: \Bbb C^d\to\Cal B(H)$ satisfies
$$
\align
C(z)C(w)+C(w)C(z)&=0,\\
C(w)^*C(z)+C(z)C(w)^*&=\<z,w\>\bold 1  \tag{1.5}
\endalign
$$
for $z,w\in\Bbb C^d$, $\<z,w\>$ denoting the 
Hermitian inner product 
$$
\<z,w\> = z_1\bar w_1+\dots+z_d\bar w_d.  
$$
The $*$-algebra generated by 
the operators $C(z)$ contains the identity and is isomorphic 
to the matrix algebra $M_{2^d}(\Bbb C)$.  

Any two irreducible representations 
of the CAR algebra (in either of its presentations 
(1.4) or (1.5)) are unitarily equivalent.  The standard 
irreducible representation of the CAR algebra is defined 
as follows.  Let $Z$ be a complex Hilbert space of 
finite dimension $d$, and let $\Lambda Z$ be the exterior 
algebra over $Z$,
$$
\Lambda Z = \Lambda^0 Z\oplus \Lambda^1 Z\oplus \Lambda^2 Z
\oplus\dots\oplus \Lambda^d Z
$$
where $\Lambda^k Z$ denotes the $k$th exterior power of $Z$.  
By definition, $\Lambda^0 Z=\Bbb C$, and the last summand 
$\Lambda^d Z$ is also isomorphic to $\Bbb C$.  
$\Lambda^k Z$ is spanned by vectors of the form 
$z_1\wedge z_2\wedge\dots\wedge z_k$, $z_k\in Z$, 
and the natural inner 
product on $\Lambda^k Z$ satisfies 
$$
\<z_1\wedge\dots\wedge z_k,w_1\wedge\dots\wedge w_k\> 
= \det(\<z_i,w_j\>),
$$
the right side denoting the determinant of the $k\times k$ 
matrix of inner products $a_{ij}=\<z_i,w_j\>$.  
$\Lambda Z$ is a direct sum of the (complex) Hilbert spaces
$\Lambda^k Z$, and it is a Hilbert space 
of complex dimension $2^d$.  

For $z\in Z$, the creation operator $C(z)$ maps $\Lambda^k Z$ 
to $\Lambda^{k+1}Z$, and acts on the generators 
as follows
$$
C(z): x_1\wedge\dots\wedge x_k\mapsto z\wedge x_1\wedge\dots\wedge x_k.  
$$
$C: Z\to\Cal B(\Lambda Z)$ is an irreducible representation 
of the canonical anticommutation relations (1.5).  
One obtains the standard irreducible Clifford 
structure (1.1) by taking the real part of this 
representation $R(z)=C(z)+C(z)^*$.  

\remark{Remarks}
In the next section we will define Dirac operators in terms of 
Clifford structures.  Because of the correspondence cited above, 
we could just as well have formulated this notion in 
terms of the anticommutation relations, avoiding Clifford structures 
entirely.  We have chosen to use them because Clifford algebras 
are associated with the Dirac operators of 
Riemannian geometry (and perhaps also for reasons of taste, 
the single equation 
(1.1) being twice as elegant as the two equations of (1.5)).  
On the other hand, we have found that proofs seem to go more 
smoothly with the anticommutation relations (1.5).  
The preceding 
observations show that nothing is lost in passing 
back and forth as needed.  
\endremark

We also want to emphasize that with any representation of 
either the Clifford relations (1.1) or the anticommutation 
relations (1.5) on a Hilbert space there are additional 
objects that are naturally associated with them, 
namely a gauge group, a number operator, and 
a $\Bbb Z_2$-grading of $H$.  By a $\Bbb Z_2$-grading of a 
Hilbert space $H$ we  simply mean a decomposition 
of $H$ into two mutually orthogonal subspaces 
$$
H = H_+\oplus H_-.  
$$
Vectors in $H_+$ (resp. $H_-$) are called even (resp. odd).  
An operator  $A\in\Cal B(H)$ is said to be of 
odd degree if $AH_+\subseteq H_-$ 
and $AH_-\subseteq H_+$, and the set of all such $A$ is a 
self-adjoint linear subspace of $\Cal B(H)$.  

\proclaim{Proposition A}Let $R: \Bbb C^d\to\Cal B(H)$ be a 
Clifford structure (1.1), and let $\Cal A$ be the 
finite dimensional $C^*$-algebra generated by the range of $R$.  
There is a unique strongly continuous unitary representation $\Gamma$ 
of the circle group $\Bbb T$ on $H$ satisfying 
$$
\align
\Gamma(\Bbb T)&\subseteq \Cal A\\
\Gamma(\lambda)R(z)\Gamma(\lambda)^*&=R(\lambda z), \qquad 
\lambda\in\Bbb T, \quad z\in \Bbb C^d, 
\endalign
$$
and such that the spectrum of $\Gamma$ starts at $0$ in 
the sense that the spectral subspaces 
$$
H_n=\{\xi\in H:  \Gamma(\lambda)\xi=\lambda^n\xi
\text{ for all }\lambda\in\Bbb T\}, \qquad n\in\Bbb Z
$$
satisfy $H_n=\{0\}$ for negative $n$ and $H_0\neq \{0\}$.  

The number operator $N$ is defined as the generator of the gauge 
group
$$
\Gamma(e^{it}) = e^{itN}, \qquad t\in\Bbb R,
$$
and is a self adjoint element of $\Cal A$ having spectrum 
$\{0,1,2,\dots,d\}$.  The $\Bbb Z_2$-grading of $H$ is 
defined by 
$$
H_+=\sum_{n {\text{ even}}} H_n, \qquad 
H_-=\sum_{n {\text{ odd}}} H_n.  
$$
\endproclaim

\demo{proof}
This is a reformulation of standard results that are 
perhaps most familiar when formulated in terms of the 
anticommutation relations.  
One may check the validity of the proposition 
explicitly for the irreducible representation 
on $\Lambda \Bbb C^d$ described above.  Since every 
Clifford structure is unitarily equivalent to a direct sum of 
copies of this irreducible one, Proposition A persists in 
the general case.  
\qedd
\enddemo

\remark{Remark 1.6}
One can single out these objects most explicitly 
in terms of the anticommutation relations $C: Z\to\Cal B(H)$ 
(1.5) over any $d$-dimensional one-particle space $Z$.  
Here, $\Cal A$ is the $C^*$-algebra 
generated by $C(Z)$ and $\Gamma$ should satisfy 
$\Gamma(\lambda)C(z)\Gamma(\lambda)^*=\lambda C(z)$
for $z\in Z$, $\lambda\in\Bbb T$, along with the two 
requirements that 
(1) the spectrum of $\Gamma$ should start at $0$ 
and (2) the gauge automorphisms of $\Cal B(H)$ should be 
inner in the sense that $\Gamma(\Bbb T)\subseteq \Cal A$.  
The number operator and 
gauge group are given by 
$$
N = C(e_1)C(e_1)^*+\dots+C(e_d)C(e_d)^*, \qquad
\Gamma(e^{it})=e^{itN}, \qquad t\in\Bbb R
$$ 
$e_1,\dots,e_d$ being any orthonormal basis for the complex Hilbert 
space $Z$.  The $\Bbb Z_2$ grading is defined by the spectral
subspaces of $\Gamma$ (or equivalently, of $N$) as in Proposition A.  
\endremark

\subhead{2.  Dirac operators and Taylor invertibility}
\endsubhead

A Dirac operator is a self-adjoint operator $D$ acting on 
a Hilbert space $H$ that has been endowed with a distinguished 
Clifford structure (1.1), satisfying three 
additional conditions.  In order to keep the bookkeeping explicit, 
we include the Clifford structure as part of the definition.  

\proclaim{Definition}
A Dirac operator of dimension $d$ is 
a pair $(D,R)$ consisting of a bounded 
self-adjoint operator $D$ acting on a Hilbert space $H$ and 
a Clifford structure $R:\Bbb C^d\to\Cal B(H)$, satisfying 
\roster
\item"{(D1)}"
(symmetry about $0$): $\Gamma(-1)D\Gamma(-1)^*=-D$, 
\item"{(D2)}"
(invariance of the Laplacian): $\Gamma(\lambda)D^2\Gamma(\lambda)^*=D^2$,
 $\lambda\in\Bbb T$, 
\item"{(D3)}"
$R(z)D + DR(z)\in \Cal A^\prime, \quad z\in\Bbb C^d$, 
\endroster
where $\Gamma: \Bbb T\to\Cal B(H)$ 
is the gauge group associated 
with $R$, and $\Cal A$ is the 
$C^*$-algebra generated by the range of $R$.  
\endproclaim

\remark{Remarks}
Let $H=H_+\oplus H_-$ be the $\Bbb Z_2$-grading of $H$
induced by the gauge group.  
(D1) is equivalent to requiring that 
$D H_+\subseteq H_-$ and $D H_-\subseteq H_+$, 
i.e., that $D$ should be an operator of odd degree.  
(D2) implies that the ``Laplacian" $D^2$ associated with $D$ should 
be invariant under the action of the gauge group 
as automorphisms of $\Cal B(H)$.  
(D3) asserts that the ``partial derivatives" of $D$ must 
commute with the operators in $R(\Bbb C^d)$.  

We have already pointed out that Clifford structures are interchangeable 
with representations $C$ of the anticommutation relations (1.5).  In terms 
of $C$, the definition of Dirac operator would be similar except that 
(D3) would be replaced with the following: 
$C(z)D+DC(z)\in\Cal A^\prime$, for every 
$z\in\Bbb C^d$.  
\endremark

There is a natural notion of isomorphism for Dirac operators, 
namely $(D,R)$ (acting on $H$) is isomorphic to $(D^\prime,R^\prime)$
(acting on $H^\prime$) if there is a unitary operator 
$U: H\to H^\prime$ such that $UD=D^\prime U$ and 
$UR(z)=R^\prime(z)U$ for every $z\in\Bbb C^d$.  Notice that 
the spectrum and multiplicity function of a Dirac operator 
are invariant under isomorphism, but of course the 
notion of isomorphism 
involves more than simple unitary equivalence of the 
operators $D$ and $D^\prime$.  

We first show how to construct a Dirac operator, starting 
with a multioperator $(T_1,\dots,T_d)$.  
Let $T_1,\dots,T_d\in\Cal B(H)$ 
be a commuting $d$-tuple of bounded operators, let 
$Z$ be a $d$-dimensional Hilbert space (which may be 
thought of as $\Bbb C^d$), and 
let $C_0: Z\to \Lambda Z$ be the irreducible representation 
of the anticommutation relations (1.5) that was described in 
section 1.  

Consider the Hilbert space $\tilde H = H\otimes\Lambda Z$
and let $C(z) = \bold 1_H\otimes C_0(z)$, $z\in Z$.  $C$ 
obviously satisfies (1.5).  Fix any orthonormal basis 
$e_1,\dots,e_d$ for $Z$ and define an operator $B$ 
on $\tilde H$ as follows 
$$
B = T_1\otimes C_0(e_1) + \dots + T_d\otimes C_0(e_d).  
$$
The pair $(D,R)$ is defined as follows
$$
D = B + B^*, \qquad R(z) = C(z)+C(z)^*, \qquad z\in Z.  
\tag{2.1}
$$
If we use the orthonormal basis to identify $Z$ with $\Bbb C^d$, 
the discussion of section 1 shows that $R$ satisfies (1.1).  

\proclaim{Proposition}
$(D,R)$ is a Dirac operator on $\tilde H$.  For 
$\lambda=(\lambda_1,\dots,\lambda_d)\in\Bbb C^d$, the Dirac 
operator of the translated $d$-tuple 
$(T_1-\lambda_1\bold 1,\dots,T_d-\lambda_d\bold 1)$ is 
$(D_\lambda,R)$, where $D_\lambda=D-R(\lambda)$.  
\endproclaim

\demo{proof}
Noting that the gauge group $\Gamma$ is related to 
$B$ by way of 
$$
\Gamma(\lambda)B\Gamma(\lambda)^*=\lambda B,  
\qquad \lambda\in\Bbb T, 
\tag{2.2}
$$
we find that 
$$
\Gamma(\lambda)D\Gamma(\lambda)^*=\lambda B + \bar\lambda B^*
$$
from which (D1) follows.  
(D3) follows after a straightforward computation 
using the anticommutation relations (1.4).  In order to check (D2), 
notice first that $B^2=0$.  Indeed, one has 
$$
B^2 = \sum_{i,j=1}^d T_iT_j\otimes C_0(e_i)C_0(e_j).  
$$
Since $T_iT_j = T_jT_i$ whereas $C_0(e_i)C_0(e_j)=-C_0(e_j)C_0(e_i)$, 
this sum must vanish.  

It follows that $D^2 = B^*B+BB^*$.  By (2.2), both $BB^*$ and 
$B^*B$ commute with the gauge group, hence so does $D^2$.  The last
sentence is immediate from (2.1).  
\qedd\enddemo

\remark{Remark}
A routine verification shows that the isomorphism 
class of this Dirac operator $(D,R)$ does not depend on 
the choice of orthonormal basis, and depends only on 
the commuting $d$-tuple $\bar T=(T_1,\dots,T_d)$.  For this 
reason we sometimes write $D_{\bar T}$ rather than $(D,R)$,  
for the Dirac operator constructed from a multioperator 
$\bar T$.  
\endremark

\remark{Comments on homology, cohomology and the Taylor spectrum}
Joseph Taylor \cite{14} introduced a notion of invertibility (and therefore
joint spectrum) for commuting $d$-tuples of operators 
$T_1,\dots,T_d$ acting on a complex Banach space.  
Taylor's notion of invertibility can be 
formulated as follows.  Let 
$$
\tilde H=\tilde H_0\oplus \tilde H_1\oplus\dots\oplus \tilde H_d
$$
be the natural decomposition of $\tilde H=H\otimes \Lambda Z$ induced by 
the decomposition of the exterior algebra $\Lambda Z$ into homogeneous forms 
of degree $k=0, 1,\dots,d$
$$
\tilde H_k = H\otimes \Lambda^k Z.  
$$
The operator 
$B=T_1\otimes c_1+\dots+T_d\otimes c_d$ of formula (1.6) satisfies 
$$
B\tilde H_k\subseteq \tilde H_{k+1}
$$
and as we have already pointed out, $B^2=0$.  Thus, the pair 
$\tilde H, B$ defines a complex (the Koszul complex of 
the $\Bbb C[z_1,\dots,z_d]$-module $H$), and when the range of $B$ is 
closed and of finite codimension in $\ker B$, we can define 
the cohomology of this complex.  Taylor defines the underlying 
$d$-tuple to be {\it invertible} if the cohomology is trivial:
$B\tilde H = \ker B$.  As we will see presently,   
for Hilbert spaces invertibility 
becomes a concrete property of the Dirac operator: 
{\it a $d$-tuple of commuting operators on $H$ is 
Taylor-invertible if and only if its Dirac operator $D$ is invertible
in $\Cal B(H\otimes\Lambda\Bbb C^d)$}.  

The Taylor spectrum of a commuting $d$-tuple $\bar T=(T_1,\dots,T_d)$ 
is defined as the set of all complex $d$-tuples 
$\lambda=(\lambda_1,\dots,\lambda_d)\in \Bbb C^d$ with the 
property that the translated $d$-tuple 
$$
(T_1-\lambda_1\bold 1,\dots,T_d-\lambda_d\bold 1)
$$
is not invertible.  In terms of the Dirac operator $(D,R)$ of 
$\bar T$, this is the set of all $\lambda\in\Bbb C^d$ such that 
$D-R(\lambda)$ is not invertible.  The relation between 
this ``Clifford spectrum" and the ordinary spectrum of $D$ 
is not very well understood.

The Taylor spectrum and Taylor's notion of invertibility are important
not only because they lead to the ``right" theorems about the spectrum 
in multivariable operator theory (see \cite{6}), 
but also and perhaps more significantly, 
because they embody the correct multivariable generalization of 
classical spectral theory as it is defined in terms of solving 
linear equations.  

In order to discuss the latter it is necessary to cast Taylor's 
cohomological picture of the joint spectrum into 
a homological picture; once that is done, 
a clear interpretation of the Taylor spectrum will emerge in terms 
of solving linear equations.  In more detail, consider the 
canonical anticommutation relations in the form (1.4) and let 
$c_1,\dots,c_d$ be the irreducible representation described in 
section 1, where $c_i$ acts as follows on the generators of 
$\Lambda^k\Bbb C^d$ 
$$
c_i:z_1\wedge\dots\wedge z_k\mapsto e_i\wedge z_1\wedge\dots\wedge z_k, 
$$
$e_1,\dots,e_d$ denoting an orthonormal basis for $\Bbb C^d$.  
Starting with a commuting $d$-tuple $T_1,\dots,T_d\in\Cal B(H)$, we 
have defined a cohomological boundary operator on $H\otimes\Lambda\Bbb C^d$
by
$$
B = T_1\otimes c_1+\dots+ T_d\otimes c_d.  
$$
Instead, let us consider the homological boundary operator 
$$
\tilde B = T_1\otimes c_1^*+\dots+ T_d\otimes c_d^*.  \tag{2.3}
$$
Formula (2.1) defines a Dirac operator $(D,R)$, and we now show
that the operators 
$$
\tilde D = \tilde B + \tilde B^*, \qquad \tilde R(z)=R(\bar z), 
\qquad z\in\Bbb C^d
$$
also define a Dirac operator $(\tilde D,\tilde R)$, 
$R$ being the Clifford structure of (2.1) and $\bar z$ denoting 
the natural conjugation in $\Bbb C^d$, 
for $z=(z_1,\dots,z_d)$, $\bar z=(\bar z_1,\dots,\bar z_d)$. 

\proclaim{Proposition: homology vs. cohomology}
The pair $(\tilde D,\tilde R)$ is a Dirac operator on 
$H\otimes\Lambda\Bbb C^d$, and it is isomorphic to the Dirac operator 
$(D,R)$ of (2.1).  The gauge group $\tilde\Gamma$ of $(\tilde D,\tilde R)$
is related to the gauge group $\Gamma$ of $(D,R)$ by 
$\tilde\Gamma(\lambda)=\lambda^d\Gamma(\lambda^{-1})$.  
\endproclaim

\demo{proof}
Consider the 
annihilation opertors $a_k=c_k^*$, $1\leq k\leq d$.  Obviously, the 
operators $a_1,\dots,a_d$ and their adjoints 
form an irreducible set of operators 
satisfying (1.4), hence there is a unitary operator $U\in\Cal B(\Lambda\Bbb C^d)$
such that $Uc_kU^*=c_k^*$, $k=1,\dots,d$.  Letting $C_0$ and $\tilde C_0$ 
be the corresponding anticommutation relations in the form (1.5), 
$$
C_0(z) = z_1c_1+\dots+z_dc_d, \qquad \tilde C_0(z) = z_1c_1^*+\dots+z_dc_d^*,
$$
we have $\tilde C_0(z)=C_0(\bar z)^*$,  and moreover 
$\tilde C_0(z)=UC_0(z)U^*$, $z\in\Bbb C^d$.  
It follows that the unitary operator 
$W=\bold 1_H\otimes U\in\Cal B(H\otimes\Lambda\Bbb C^d)$ 
satisfies $WC(z)W^*=C(\bar z)^*$, 
$z\in\Bbb C^d$.  Since 
$$
\tilde R(z)=R(\bar z)=C(\bar z)+C(\bar z)^* =
W(C(z)^*+C(z))W^*=WR(z)W^*
$$
and since $\tilde B=WBW^*$, $W$ implements an
isomorphism of the pair $(D,R)$ and the pair $(\tilde D,\tilde R)$.  
Thus, $(\tilde D,\tilde R)$ is a 
Dirac operator isomorphic to $(D,R)$.

Letting $C_k=\bold 1\otimes c_k$, $k=1,\dots,d$ the 
number operators $\tilde N$ and $N$ of $(\tilde D,\tilde R)$ and 
$(D,R)$ are seen to be 
$$
\tilde N=C_1^*C_1+\dots+C_d^*C_d,\quad N=C_1C_1^*+\dots+C_dC_d^*,
$$
so by the anticommutation relations (1.4) we have 
$\tilde N=d\cdot \bold 1-N$, and the formula relating 
$\tilde \Gamma$ to $\Gamma$ follows from Remark 1.6.  
\qedd\enddemo

In particular, the preceding proposition implies that the 
Taylor spectrum can be defined in either cohomological 
terms (using $(D,R)$ and its associated coboundary operator $B$) 
or in homological terms (using $(\tilde D,\tilde R)$ and its 
boundary operator $\tilde B$).  It 
is the homological formulation that leads to the following 
interpretation.  

Classical spectral theory starts with 
the problem of solving linear equations
of the form $Tx = y$, where $T$ is a given operator in $\Cal B(H)$, 
$y$ is a given vector in $H$, and $x$ is to be found; $T$ is said 
to be invertible when for every $y$ there is a unique $x$.  
Taylor's notion of invertibility in its {\it homological} form 
provides the 
correct generalization to higher dimensions 
of this fundamental notion in dimension one.   
In dimension two for example, one has a pair $T_1,T_2$ 
of commuting operators acting on a Hilbert space $H$, and one is interested 
in solving equations of the form 
$$
T_1x_1 + T_2x_2 = y, \tag{2.4}
$$
where $y$ is a given vector in $H$.  Of course the pair $(x_1,x_2)$ 
is never uniquely determined by $y$, since if $(x_1,x_2)$ solves this 
equation then so does $(x_1^\prime, x_2^\prime)$ where 
$x_1^\prime = x_1+T_2\zeta$ and $x_2^\prime=x_2-T_1\zeta$ where 
$\zeta\in H$ is arbitrary.  Equivalently, 
$$
\align
x_1^\prime &=x_1 + T_1\xi_{11}+T_2\xi_{12}\\
x_2^\prime &=x_2 + T_1\xi_{21}+T_2\xi_{22}, 
\endalign
$$
where the vectors $\xi_{ij}$, $1\leq i,j\leq 2$ satisfy 
$\xi_{ji}=-\xi_{ij}$ for all $i,j$ but are otherwise arbitrary (of 
course, we must have $\xi_{11}=\xi_{22}=0$ and $\xi_{12}=\zeta$).  
Such perturbations $(x_1^\prime,x_2^\prime)$ can be 
written down 
independently of any properties of the given operators $T_1$, $T_2$
(beyond commutativity, of course), and for that reason we 
will call them  {\it tautological} perturbations
of the given solution $x_1,x_2$.  
In order to understand 
how to solve such equations one needs to determine 
what happens {\it modulo} tautological perturbations, and 
for that one must look at the homology of (2.4).  

Since we are in dimension two we can write 
$$
H\otimes\Lambda\Bbb C^2= \Omega_0\oplus\Omega_1\oplus \Omega_2
$$
where $\Omega_0=H$, $\Omega_1=\{(x_1,x_2): x_k\in H\}$, and 
$\Omega_2$ is parameterized as a space of ``antisymmetric" sequences 
as follows
$$
\Omega_2=\{(\xi_{ij}): 1\leq i,j\leq 2, \xi_{ij}=-\xi_{ji} 
{\text{ for all }}i,j\}.  
$$
Of course, $\Omega_2$ is isomorphic to $H$ by way of the map which 
associates to a vector $\zeta\in H$ the antisymmetric sequence 
$\xi_{11}=\xi_{22}=0, \xi_{12}=\zeta, \xi_{21}=-\zeta$.  
The homological boundary operator $B=T_1\otimes c_1^*+T_2\otimes c_2^*$ 
of the complex 
$$
0\longleftarrow\Omega_0\longleftarrow\Omega_1\longleftarrow\Omega_2
\longleftarrow 0
\tag{2.5}
$$
acts as follows.  On $\Omega_1$, 
$B(x_1,x_2)=T_1x_1+T_2x_2$, and on $\Omega_2$
$$
B(\xi_{ij})=(T_1\xi_{11}+T_2\xi_{12},T_1\xi_{21}+T_2\xi_{22}) =
(T_2\xi_{12},-T_1\xi_{12}).
$$  
Apparently, (2.4) has a solution iff $y$ belongs 
to $B\Omega_1=T_1H+T_2H$.  Given a solution $(x_1,x_2)$ of (2.4) and another 
pair of vectors $(x_1^\prime,x_2^\prime)$, $(x_1^\prime,x_2^\prime)$ is 
also a solution iff the difference $(x_1-x_1^\prime,x_2-x_2^\prime)$ 
belongs to $\ker B$.  Given that $(x_1^\prime,x_2^\prime)$ is 
a solution, then it is a tautological perturbation of $(x_1,x_2)$ iff 
the difference $(x_1-x_1^\prime,x_2-x_2^\prime)$  belongs to $B\Omega_2$.  
Finally, the kernel of the boundary operator at $\Omega_2$ is 
identified with $\ker T_1\cap\ker T_2$.  We conclude that the 
complex (2.5) is exact iff a) $T_1H+T_2H=H$, b) $\ker T_1\cap\ker T_2=\{0\}$, 
and c) solutions of (2.4) are unique up to tautological perturbations.  
While the algebra is more subtle in higher dimensions the fundamental issues 
are the same, and that is why the Taylor spectrum is important in 
multivariable spectral theory.  

We will not have to delve into homological issues here; 
but the above comments do show
that the theory of abstract Dirac operators is rooted in concrete 
problems of linear algebra that are associated with solving linear 
equations involving commuting sets of operators. 
\endremark

Taylor's definition of invertibility can be reformulated 
in terms of the Dirac operator $D_{\bar T}$, and then extended 
to define Fredholm $d$-tuples and their index.  
In more detail, in the proof of the previous proposition 
we have already pointed out that $D^2 = B^*B + BB^*$; and since 
$B\tilde H$ and $B^*\tilde H$ are orthogonal, we conclude that 
$B\tilde H =\ker B$ iff $D^2$ is invertible.  

Conclusion: A commuting $d$-tuple $(T_1,\dots, T_d)$ is invertible 
if and only if its Dirac operator is invertible.  

By a Fredholm $d$-tuple we mean one whose Dirac operator $(D,R)$ is 
Fredholm in the sense that the self-adjoint operator $D$
has closed range and finite dimensional 
kernel.  The index of a Fredholm $d$-tuple is defined as 
follows.  By property (D1) we have 
$D\tilde H_+\subseteq \tilde H_-$ and 
 $D\tilde H_-\subseteq \tilde H_+$.  Thus we may consider the operator 
$$
D_+ = D\restriction_{H_+}\in \Cal B(\tilde H_+, \tilde H_-),  
$$
whose adjoint is given by 
$$
D_+^* = D\restriction_{H_-}\in \Cal B(\tilde H_-, \tilde H_+).  
$$
For a Fredholm $d$-tuple $\bar T=(T_1,\dots,T_d)$, 
$D_+$ is a Fredholm operator from 
$\tilde H_+$ to $\tilde H_-$, and the
index of $\bar T$ is defined by 
$$
{\text{ind}}(\bar T)=\dim\ker(D_+) - \dim\ker(D_+^*).
$$  
One can define semi-Fredholm $d$-tuples similarly, but we 
do not require the generalization here.

\subhead{3.  Dirac operators and Hilbert modules over $\Bbb C[z_1,\dots,z_d]$}
\endsubhead

In this section we prove the following result, which implies 
that Dirac operators $(D,R)$ contain exactly the same geometric 
information as multioperators $\bar T$.  

\proclaim{Theorem A}
For every $d$-dimensional Dirac operator $(D,R)$ 
there is a commuting $d$-tuple 
$\bar T=(T_1,\dots,T_d)$ acting on a Hilbert space $H$ 
such that $(D,R)$ is isomorphic to $D_{\bar T}$.  
If $\bar T^\prime=(T_1^\prime,\dots,T_d^\prime)$ is another 
commuting $d$-tuple acting on $H^\prime$, then $D_{\bar T}$ and 
$D_{\bar T^\prime}$ are isomorphic if and only if 
there is a unitary operator $U: H\to H^\prime$ such that 
$UT_k=T_k^\prime U$ for every $k=1,\dots,d$.  
\endproclaim

\demo{proof}
Let $K$ be the underlying Hilbert space of 
$(D,R)$, so that $D=D^*\in\Cal B(K)$ and 
$R: \Bbb C^d\to\Cal B(K)$ is a 
Clifford structure (1.1) which satisfy 
(D1), (D2), (D3).  

Consider the map $C: \Bbb C^d\to\Cal B(K)$ defined by 
$C(z)=(1/2)(R(z)-iR(iz))$.  The discussion of section 1 
implies that $C$ satisfies the 
anticommutation relations (1.5), 
and  $R(z)=C(z) + C(z)^*$.    
$C$ is unitarily equivalent to 
a direct sum of copies of the standard irreducible representation 
$C_0$ of the anticommutation relations on $\Lambda \Bbb C^d$;
thus by replacing $(D,R)$ with an isomorphic copy we may 
assume that there is a Hilbert space $H$ such that 
$K=H\otimes \Lambda\Bbb C^d$ and that 
$R(z) = C(z)+C(z)^*$ where $C(z)$ is defined on 
$H\otimes\Lambda \Bbb C^d$ by
$$
C(z) = \bold 1_H\otimes C_0(z), \qquad z\in\Bbb C^d.  
$$
We must exhibit a {\it commuting} set of operators $T_1,\dots,T_d$ 
on $H$ so that $D=B+B^*$ where 
$$
B=T_1\otimes C_0(e_1)+\dots+ T_d\otimes C_0(e_d), 
$$
$e_1,\dots,e_d$ being the usual orthonormal basis for $\Bbb C^d$.  

To that end, let $\Cal A$ be the finite dimensional $C^*$-algebra 
$\Cal A=\bold 1_H\otimes\Cal B(\Lambda\Bbb C^d)$.  
The $C^*$-algebras generated by $R(\Bbb C^d)$ and 
$C(\Bbb C^d)$ are the same, and in fact 
$$
C^*(R(\Bbb C^d))=C^*(C(\Bbb C^d))=\Cal A.   \tag{3.1}.  
$$
By (D1), 
$R(z)D+DR(z)$ must commute with 
$\Cal A$ for every $z\in \Bbb C^d$ and, 
in view of the relation $C(e_k)^* = 2(R(e_k)+iR(ie_k))$ we have 
$C(e_k)^*D+D C(e_k)^*\in\Cal A^\prime$.  Thus for every $k$ there is 
a unique operator $T_k\in\Cal B(H)$ such that 
$$
C(e_k)^*D + DC(e_k)^* = T_k\otimes\bold 1_{\Lambda\Bbb C^d}.  \tag{3.2}
$$

For each $k=1,\dots,d$, let $c_k=C_0(e_k)\in\Cal B(\Lambda\Bbb C^d)$, 
and consider the operator
$$
B=T_1\otimes c_1+\dots +T_d\otimes c_d\in\Cal B(H\otimes\Lambda\Bbb C^d).
\tag{3.3}
$$
In order to show that $D=B+B^*$ we will make use of 

\proclaim{Lemma}
Let $R: \Bbb C^d\to\Cal B(K)$ be a 
Clifford structure (1.1) on $K$ and let 
$\Gamma:\Bbb T\to\Cal B(K)$ 
be its gauge group.  
Every operator $A\in\Cal B(K)$ 
satisfying $R(z)A+AR(z)=0$ for every $z\in\Bbb C^d$
admits a decomposition $A = A_0\Gamma(-1)$, where $A_0$ belongs to 
the commutant of $C^*(R(\Bbb C^d))$.   
In particular, such an operator must also be gauge invariant 
in the sense that $\Gamma(\lambda)A\Gamma(\lambda)^*=A$, 
$\lambda\in\Bbb T$.  
\endproclaim

\demo{proof}
Since $\Gamma(-1)R(z)\Gamma(-1)^*=R(-z)=-R(z)$ it follows 
that $\Gamma(-1)$ anticommutes with $R(z)$ for every 
$z\in\Bbb C^d$.  Since $A$ also anticommutes with $R(z)$, 
the operator $A_0=A\Gamma(-1)$ must 
commute with $R(z)$, and we have 
$A=A\Gamma(-1)^2=A_0\Gamma(-1)$ as required. 
The last assertion follows
from this decomposition, because for every $\lambda\in\Bbb T$, 
$\Gamma(\lambda)$ belongs to
the $C^*$-algebra generated by the range of $R$ and hence 
commutes with both factors $A_0$ and $\Gamma(-1)$.  
\qedd\enddemo

We now show that for $B$ as in (3.3) we have $D=B+B^*$.  
Indeed, since $C(e_k)=\bold 1_H\otimes c_k$ and the 
$c_k$ satisfy the anticommutation relations (1.4) we have 
$$
\align
C(e_k)B + BC(e_k) &= \sum_{j=1}^d T_j\otimes (c_kc_j+c_jc_k)=0,\\
C(e_k)^*B + BC(e_k)^* &= \sum_{j=1}^d T_j\otimes (c_k^*c_j+c_jc_k^*)
=\sum_{j=1}^d T_j\otimes \delta_{jk}\bold 1=T_k\otimes\bold 1.  
\endalign
$$
Using the definition of $T_k$ (3.2) it follows from the preceding 
calculation that the difference $D-B-B^*$ must anticommute with all of the 
operators $C(e_j), C(e_k)^*$, $1\leq j,k\leq d$.  Since 
$R(z) = C(z)+C(z)^*$ it follows that $D-B-B^*$ anticommutes 
with $R(z)$ for every $z\in \Bbb C^d$.  

By the Lemma, there is a (necessarily unique) 
operator $X\in\Cal B(H)$ such that 
$$
D-B-B^* = X\otimes\Gamma_0(-1), \tag{3.4}
$$
where $\Gamma_0:\Bbb T\to\Cal B(\Lambda\Bbb C^d)$ 
is the natural gauge action on $\Lambda\Bbb C^d$ and 
$\Gamma(\lambda) = \bold 1_H\otimes\Gamma_0(\lambda)$.  
We want to show that $X=0$.  For that, recall that $D$ is 
odd (property (D1)) and $B$ is clearly odd by its 
definition (3.3).  Hence $D-B-B^*$ is odd, so it must anticommute 
with the unitary operator 
$\Gamma(-1)=P_{H_+}-P_{H_-}$.  On the other hand (3.4) 
implies that it commutes with $\Gamma(-1)$.  Since 
$\Gamma(-1)$ is invertible, $D-B-B^*=0$.

What remains to be proved is that the operators $T_k$ of (3.2) 
commute with each other.  Indeed, we claim first 
that $B^2=0$.  Since we have established that 
$D=B+B^*$ we can write 
$$
D^2 = B^*B + BB^* + B^2 + {B^*}^2.  \tag{3.5}
$$
From the definition of $B$ (3.3) we have 
$$
\Gamma(\lambda)B\Gamma(\lambda)^* = 
\sum_{k=1}^d T_k\otimes \Gamma_0(\lambda)c_k\Gamma_0(\lambda)^* =
\lambda\sum_{k=1}^d T_k\otimes c_k=\lambda B,\qquad \lambda\in\Bbb T.  
$$
It follows that $B^*B$ and $BB^*$ are invariant under the 
action of the gauge group, and that 
$\Gamma(\lambda)B^2\Gamma(\lambda)^*=\lambda^2B^2$.  Thus
$$
\Gamma(\lambda)D^2\Gamma(\lambda)^* = 
B^*B + BB^* + \lambda^2 B^2 + \bar\lambda ^2{B^*}^2.  \tag{3.6}
$$
Because of (D2), 
the left side of (3.6) does not depend on $\lambda$.  Hence 
by equating Fourier coefficients on left and right we find 
that $B^2 = {B^*}^2 = 0$.  

We can now show that the operators $T_k$ defined by (3.2) mutually 
commute.  Consider the operator 
$C$ defined on $H\otimes \Lambda\Bbb C^d$ by 
$$
C=\sum_{1\leq j<k\leq d}(T_jT_k-T_kT_j)\otimes c_jc_k.  \tag{3.7}
$$
Since the operators 
$\{c_jc_k: 1\leq j<k\leq d\}\subseteq\Cal B(\Lambda\Bbb C^d)$ 
are linearly independent, it is enough to show that 
$C=0$.  To see this, we use the anticommutation 
relations $c_kc_j+c_jc_k=\delta_{jk}\bold 1$ to write 
$$
\align
C &= \sum_{1\leq j<k\leq d}T_jT_k\otimes c_jc_k - 
\sum_{1\leq j<k\leq d}T_kT_j\otimes c_jc_k \\
&=\sum_{1\leq j<k\leq d}T_jT_k\otimes c_jc_k +
\sum_{1\leq j<k\leq d}T_kT_j\otimes c_kc_j\\
&=\sum_{1\leq p,q\leq d}T_pT_q\otimes c_pc_q = B^2=0.  
\endalign
$$
That completes the proof that every Dirac operator 
is associated with a commuting $d$-tuple.  

Suppose now that we are given two commuting 
$d$-tuples $\bar T$ and $\bar T^\prime$, acting on 
Hilbert spaces $H$ and $H^\prime$.   It is obvious that 
if $U:H\to H^\prime$ is a unitary operator satisfying 
$UT_k=T_k^\prime U$ for every $k=1,\dots,d$, then 
$W=U\otimes\bold 1: H\otimes\Lambda\Bbb C^d\to 
H^\prime\otimes\Lambda\Bbb C^d$ is a unitary operator 
which implements an isomorphism of the respective 
Dirac operators.

Conversely, let 
$W: H\otimes\Lambda\Bbb C^d\to H^\prime\otimes\Lambda\Bbb C^d$
be a unitary operator implementing an isomorphism of the 
respective Dirac operators $(D,R)$ and $(D^\prime,R^\prime)$
associated with $\bar T$ and $\bar T^\prime$.  
Let $R_0: \Bbb C^d\to\Lambda\Bbb C^d$ be the irreducible 
Clifford structure defined in section 1.  Since 
$R(z)=\bold 1_H\otimes R_0(z)$ and 
$R^\prime(z)=\bold 1_{H^\prime}\otimes R_0(z)$, it follows 
that $H$ and $H^\prime$ have the same dimension (namely 
the common multiplicity of the 
unitarily equivalent Clifford structures 
$R$ and $R^\prime=WRW^*$).  Thus by replacing $\bar T^\prime$ 
with a unitarily equivalent $d$-tuple, we can assume 
that $H=H^\prime$, i.e., that both $d$-tuples act on the same 
Hilbert space $H$.  

In these ``coordinates", the relation 
$$
W(\bold 1_H\otimes R_0(z))W^* = \bold 1_H\otimes R_0(z), 
\qquad z\in\Bbb C^d
$$
implies that $W$ commutes with 
$\bold 1_H\otimes\Cal B(\Lambda\Bbb C^d)$, 
the $C^*$-algebra generated 
by $R(\Bbb C^d)$.  Thus $W$ decomposes 
$W=U\otimes\bold 1_{\Lambda\Bbb C^d}$ where $U$ is a 
uniquely determined unitary 
operator on $H$.  
Now according to the definition of Dirac operators (2.1), 
we have $D=B+B^*$, $D^\prime=B^\prime+{B^\prime}^*$, 
where 
$$
B=T_1\otimes c_1+\dots+T_d\otimes c_d,\qquad
B^\prime=T_1^\prime\otimes c_1+\dots+T_d^\prime\otimes c_d,
$$
$c_1,\dots,c_d$ being the irrecucible representation of 
the canonical anticommutation relations (1.4) associated with 
$R_0$.  Letting $C_k=\bold 1_H\otimes c_k$ and using (1.4),
a routine calculation gives
$$
C_k^*D + DC_k^*=T_k\otimes\bold 1,
\quad 
C_k^*D^\prime + D^\prime C_k^*=T_k^\prime\otimes\bold 1,
\qquad k=1,\dots,d.  
$$
Since $U\otimes\bold 1=W$ commutes with all $C_k^*$ and 
satisfies $WDW^*=D^\prime$, it follows that for every 
$k=1,\dots,d$ we have 
$$
UT_kU^*\otimes \bold 1=W(C_k^*D+DC_k^*)W^*=C_kD^\prime+D^\prime C_k^*
=T_k^\prime\otimes\bold 1,
$$
and hence $U$ implements a unitary equivalence of 
$\bar T$ and $\bar T^\prime$.  
That completes the proof of Theorem A.\qedd\enddemo

\remark{Remark 3.8}
It is worth pointing out that the proof of Theorem A shows 
how one may go directly from a Dirac operator $(D,R)$ (acting 
on $H$) to the 
Koszul complex of its underlying $d$-tuple $\bar T$ 
(the operators $T_1,\dots,T_d$ acting on some other Hilbert space) without 
making explicit reference to $\bar T$.  Indeed, considering the spectral 
representation of the gauge group of $R$
$$
\Gamma(\lambda)=\sum_{n=-\infty}^\infty \lambda^n E_n =
\sum_{n=0}^d \lambda^n E_n,  \qquad\lambda\in\Bbb T
$$
and the operator $B=\sum_n E_{n+1}DE_n$, 
from the proof of Theorem A one finds that 
$$
B^2 = 0, \qquad D = B+B^*.  \tag{3.9}
$$
Moreover, the spectral subspaces $H_n=E_nH$ satisfy 
$BH_n\subseteq H_{n+1}$, $B^*H_n\subseteq H_{n-1}$, and the Koszul complex is
given by 
$$
0\longrightarrow H_0\longrightarrow H_1\longrightarrow 
\dots\longrightarrow H_d\longrightarrow 0
$$
with cohomology defined by $B$.  
\endremark

\subhead{4.  Stability of the Curvature invariant: graded case}
\endsubhead

Recall from \cite{3} that a commuting $d$-tuple of operators 
$(T_1,\dots,T_d)$ on a Hilbert space $H$ is said to be graded 
if it is circularly symmetric in the sense that 
there is a strongly continuous unitary representation 
$\Gamma: \Bbb T\to\Cal B(H)$ such that 
$$
\Gamma(\lambda)T_k\Gamma(\lambda)^*=\lambda T_k, \qquad 
k=1,\dots,d, \quad \lambda\in\Bbb T.  
$$
Many examples of graded $d$-contractions 
were described in 
\cite{3}; in particular, all examples of 
$d$-contractions 
that were associated with projective algebraic varieties 
(and their finitely generated modules) are graded.  

It was shown in (\cite{3}, see Theorem B) that the curvature invariant 
of a pure graded finite rank $d$-contraction is an 
integer, namely the Euler characteristic of a certain finitely 
generated algebraic module over $\Bbb C[z_1,\dots,z_d]$ that is 
associated naturally with $\tilde T$.  
However, in the ungraded case this formula fails (both sides 
of this formula still make sense in the ungraded case, but examples 
are given in \cite{3} for which they are unequal).  This led us 
to ask in \cite{3},\cite{4} if $K(\bar T)$ is an 
integer even when $\bar T$ is ungraded.  That has been 
recently proved by Greene, Richter and Sundberg \cite{9}, 
and in fact the results of \cite{9} show that 
the integer $K(\bar T)$ 
can be identified as the (almost everywhere constant) rank of the 
boundary values of a certain operator-valued 
``inner"  function that is naturally associated with 
$\bar T$ via  dilation theory (a 
fuller discussion of this inner operator and 
its relation to $\bar T$ can be found in \cite{2}).  
It is fair to say that 
the rank of this inner function is not easily 
computed in terms of the operator theory of $\bar T$, 
and thus we were led to ask if there is a formula that relates 
the curvature invariant more 
directly to the geometry of $\bar T$...preferably 
in terms of an expression that is obviously an integer.  

It is also noteworthy that the 
asymptotic formula for the curvature (Theorem C of \cite{3}) 
implies that it has certain stability properties; 
for example, the curvature is stable under the operation 
of restricting to an invariant subspace of finite codimension.  
But nothing was known about stability of the 
curvature invariant under compact perturbations.  

These questions led us to search for another 
formula for the curvature invariant that looks more like 
an index theorem in the sense that it equates the 
curvature invariant to the index of some operator.  Such 
a formula would presumably lead to stability under compact 
perturbations, it would imply that the curvature invariant 
is in all cases an integer, and it would more closely 
resemble the Gauss-Bonnet-Chern formula in its 
modern incarnation as an index theorem 
(for example, see p. 311 of \cite{8}).  As a first 
step in this direction, we offer the following.

\proclaim{Theorem B}
Let $\bar T=(T_1,\dots,T_d)$ be a pure $d$-contraction of 
finite rank acting on a Hilbert space $H$.  Assume that 
$\bar T$ is graded and let $(D,R)$ be its Dirac operator.  
Then both $\ker D_+$ and $\ker D_+^*$ are finite dimensional
and 
$$
(-1)^dK(\bar T) = \dim \ker D_+ - \dim \ker D_+^*.  
$$
\endproclaim

\remark{Remark}
Notice that we have {\it not} assumed that $D$ 
is a Fredholm operator.  However, when it is Fredholm 
we have the following stability.  
\endremark

\proclaim{Corollary}
Let $\bar T=(T_1,\dots,T_d)$  and 
$\bar T^\prime = (T_1^\prime,\dots,T_d^\prime)$ 
be two pure $d$-contractions 
of finite rank acting on respective Hilbert 
spaces $H$, $H^\prime$.   Assume that both 
$\bar T$ and $\bar T^\prime$ are graded, that 
$\bar T$ is Fredholm, and 
that they are unitarily equivalent modulo 
compacts in the sense that there is a unitary 
operator $U: H\to H^\prime$ such that 
$$
UT_k - T_k^\prime U {\text{ is compact}}, \qquad k=1,\dots,d.  
$$  
Then $K(\bar T) = K(\bar T^\prime)$.  
\endproclaim

\demo{proof of Corollary}
Let $(D,R)$ and $(D^\prime,R^\prime)$ be the Dirac operators 
of $\bar T$ and $\bar T^\prime$, acting 
on respective Hilbert spaces $\tilde H=H\otimes \Lambda\Bbb C^d$ and 
$\tilde H^\prime=H^\prime\otimes\Lambda\Bbb C^d$.  
The hypothesis implies 
that the unitary operator 
$W: U\otimes\bold 1: \tilde H\to \tilde H^\prime$
satisfies $WR(z)=\tilde R(z)W$ for all $z\in\Bbb C^d$, 
and $WD- D^\prime W$ is compact.  The first of these two 
relations implies that $W$ implements an equivalence 
of the respective gauge groups 
$W\Gamma(\lambda)=\Gamma^\prime(\lambda)W$, and hence 
$W$ carries the $\Bbb Z_2$-grading of $\tilde H$ to 
that of $\tilde H^\prime$.  It follows that the restrictions 
of $W$ to the even and odd subspaces of $\tilde H$ implement a 
unitary equivalence modulo compact operators of 
the two operators $D_+$ and $D^\prime_+$.  Since 
$D_+$ is Fredholm by hypothesis, $D_+^\prime$ must 
be Fredholm as well, and moreover they must have the 
same index.  From Theorem B we conclude that 
$K(\bar T)=K(\bar T^\prime)$.  
\qedd\enddemo

Before giving the proof of Theorem B, we recall some algebraic 
preliminaries.  Let $\Cal A$ be the complex polynomial algebra 
$\Bbb C[z_1,\dots,z_d]$.  By an $\Cal A$-module we mean a complex 
vector space $M$ which is endowed with a commuting $d$-tuple of 
linear operators $T_1,\dots,T_d$, the module structure being defined 
by $f\cdot\xi=f(T_1,\dots,T_d)\xi$, $f\in \Cal A$, $\xi\in M$.  
$M$ is said to be finitely generated if there is a finite set 
$\xi_1,\dots,\xi_s$ of vectors in $M$ such that 
$$
M = \{f_1\cdot\xi_1+\dots+f_s\cdot\xi_s: f_1,\dots,f_s\in\Cal A\}.  
$$
The free module of rank $1$ is 
defined to be $\Cal A$ itself, with the module 
action associated with multiplication of polynomials.  The free 
module of rank $r=1,2,\dots$ is the direct sum of $r$ copies of 
the free module of rank $1$, with the obvious module action on 
$r$-tuples of polynomials.  

Hilbert's Syzygy theorem implies that every finitely generated 
$\Cal A$-module has a finite free resolution \cite{7} in the 
sense that there is an exact sequence of $\Cal A$-modules 
$$
0\longrightarrow F_n\longrightarrow\dots\longrightarrow 
F_1\longrightarrow M\longrightarrow 0
\tag{4.1}
$$
where each $F_k$ is a free module of finite rank.  In \cite{3}, 
we defined the Euler characteristic of $M$ in terms of finite 
free resolutions (4.1) as follows
$$
\chi(M) = \sum_{k=1}^n (-1)^{k+1}\rank(F_k).  \tag{4.2}
$$
This integer does not depend on the particular resolution of $M$ 
chosen to define it.  

We must relate $\chi(M)$ to the alternating sum of the Betti 
numbers of the Koszul complex of $M$; since the latter is 
also called the Euler characteristic, we distinguish it 
from $\chi(M)$ by calling it the Euler number of $M$ and 
by writing it as $e(M)$.  The Euler number is defined as follows.  

The Koszul complex of an $\Cal A$-module $M$ is defined as the 
$\Cal A$-module 
$$
M\otimes\Lambda\Bbb C^d=\Omega^0\oplus\Omega^1\oplus\dots\oplus
\Omega^d, 
$$
where $\Omega^k=M\otimes\Lambda^k\Bbb C^d$ is the submodule of $k$-forms, 
with coboundary operator
$$
B=T_1\otimes c_1+\dots+T_d\otimes c_d
$$
exactly as we have done above in the case where $M$ is a Hilbert 
space and the $T_k$ are bounded linear operators.  Letting $B_k$ 
be the restriction of $B$ to $\Omega^k$ we have 
a corresponding cohomology space
$H^k(M)=\ker B_k/{\text{ran}} B_{k-1}$ for $1\leq k\leq d$, 
$H^0(M)=\ker B_0$, which may or may not be 
finite dimensional.  $M$ is said to be of {\it finite type} if 
$H^k(M)$ is finite dimensional for every $0\leq k\leq d$, and in 
that case the Euler number is defined by 
$$
e(M)=\sum_{k=0}^d(-1)^k\dim H^k(M).  \tag{4.3}
$$

Taking $M$ to be the free module $\Cal A$ of rank one, it is well-known 
that $H^k(\Cal A)=0$ for $0\leq k\leq d-1$ and that 
$H^d(\Cal A)=\Cal A/(z_1\Cal A+\dots+z_d\Cal A)\cong \Bbb C$ is one-dimensional.  
It follows that for a free 
module $F$ of arbitrary finite rank, we have 
$$
e(F)=(-1)^d\cdot \rank F.  \tag{4.4}
$$

The following result is part of the lore of commutative 
algebra; we sketch a proof for the reader's convenience.  

\proclaim{Lemma 1}
Let $0\longrightarrow K\longrightarrow L\longrightarrow M\longrightarrow 0$ 
be a short exact sequence of $\Cal A$ modules, some two of which are 
of finite type.  Then all are of finite type and we have 
$$
e(L)=e(K) + e(M).  
$$
\endproclaim

\demo{proof}
Letting $\kappa(N)$ denote the Koszul complex of an $\Cal A$-module N, 
one sees that $\kappa(N)$ has $d+1$ nonzero terms, and 
the corresponding sequence of complexes
$$
0\longrightarrow\kappa(K)\longrightarrow\kappa(L)\longrightarrow
\kappa(M)\longrightarrow 0
$$
is exact.  Thus by fundamental principles we obtain a 
long exact sequence of cohomology spaces which contains 
at most $3d + 3$ nonzero terms.  Two of any three 
consecutive terms in the latter sequence are finite dimensional because 
two of the three modules $K, L, M$ are assumed to have finite dimensional 
cohomology.  By exactness all cohomology spaces are finite dimensional 
and the alternating sum of their dimensions must be zero.  
The asserted formula follows.  
\qedd\enddemo

\proclaim{Lemma 2}
Every finitely generated $\Cal A$-module $M$ is of finite type, 
and  
$$
e(M)=(-1)^d\chi(M).  
$$
\endproclaim

\demo{proof}
Choose a finite free resolution of $M$ in the form (4.1)
$$
0\longrightarrow F_n\longrightarrow\dots\longrightarrow 
F_1\longrightarrow M\longrightarrow 0.  
$$
Let $R_k\subseteq F_{k-1}$ be the image of $F_k$ in 
$F_{k-1}$, $2\leq k\leq n$, and let $R_1\subseteq M$ 
be the image of $F_1$.  Starting at the left of 
(4.1) we have a short exact sequence of 
modules 
$$
0\longrightarrow F_n\longrightarrow 
F_{n-1}\longrightarrow R_{n-1}\longrightarrow 0, 
$$ 
the first two of which are of finite type.  By Lemma 1, 
$R_{n-1}$ is of finite type and 
$$
e(R_{n-1})=e(F_{n-1})-e(F_n).  
$$
Moving one step to the right, the same argument applied to 
$$
0\longrightarrow R_{n-1}\longrightarrow 
F_{n-2}\longrightarrow R_{n-2}\longrightarrow 0
$$
shows that $R_{n-2}$ is of finite type and 
$$
e(R_{n-2})=e(F_{n-2})-e(R_{n-1})=e(F_{n-2})-e(F_{n-1})+e(F_n).  
$$
Continuing in this way to the end of the sequence, 
we arrive at the conclusion that 
$M$ is of finite type and 
$$
e(M)=\sum_{k=1}^n(-1)^{k+1}e(F_k)=(-1)^d\sum_{k=1}^n(-1)^{k+1}\rank(F_k)
=(-1)^d\chi(M), 
$$
where in the second equality we have made use of (4.4).  
\qedd\enddemo

\demo{proof of Theorem B}
We are assuming that $\bar T$ is graded; this means that there 
is a continuous unitary representation of the circle 
group $U: \Bbb T\to\Cal B(H)$ such that 
$$
U_\lambda T_k U_\lambda^*=\lambda T_k, \qquad 1\leq k\leq d.  \tag{4.5}
$$
Let $\Delta = (\bold 1-T_1T_1^*-\dots-T_dT_d^*)^{1/2}$ be the 
defect operator of $\bar T$.  By hypothesis, $\Delta$ is of finite
rank, and the canonical algebraic module $M_H$ associated with 
$\bar T$ 
$$
M_H={\text{span}}\{f(T_1,\dots,T_d)\zeta: 
f\in\Cal A, \zeta\in\Delta H\}
$$
is a finitely generated $\Cal A$ module.  Because $\bar T$ is 
pure, $M_H$ is dense in $H$ (see \cite{3}, Proposition 5.4).  

It follows that $M_H\otimes \Lambda\Bbb C^d$ 
is dense in $\tilde H=H\otimes\Lambda\Bbb C^d$.  Let 
$D\in\Cal B(\tilde H)$ be the Dirac operator of $\bar T$.  We 
will show that both $\ker D_+$ and $\ker D_+^*$ are finite 
dimensional subspaces of $M_H\otimes\Lambda\Bbb C^d$, 
and that in fact we have 
$$
\dim\ker(D_+)=\sum_{k {\text { even}}}\dim H^k(M_H), 
\quad \dim\ker(D_+^*)=\sum_{k {\text { odd}}}\dim H^k(M_H) 
\tag{4.6}
$$
where $M_H\otimes\Lambda\Bbb C^d$ is viewed as the Koszul complex 
of $M_H$.  Assuming for the moment that (4.6) has been established 
we find that 
$$
\dim\ker D_+-\dim\ker D_+^*=\sum_{k=0}^d(-1)^k\dim H^k(M_H)=e(M_H), 
$$
and by Lemma 2 the right side is $(-1)^d\chi(M_H)$.  By Theorem 
B of \cite{3}, the latter is $(-1)^d K(H)$, and the proof 
of Theorem B above will be complete.  

In order to establish (4.6) we make use of the grading as follows.  
Let $c_1,\dots, c_d$ be operators on $\Bbb C^d$ 
satisfying the anticommutation relations (1.4) 
and let 
$$
B=T_1\otimes c_1+\dots+T_d\otimes c_d
$$
be the coboundary operator on $\tilde H$.  Since 
$D^2=B^*B + BB^*$, the kernel of $D$ is given by 
$\ker D = \ker B \cap \ker B^*$.  Let $V: \Bbb T\to\Cal B(\tilde H)$ 
be the unitary representation corresponding to $U$, 
$V_\lambda=U_\lambda\otimes 1_{\Lambda\Bbb C^d}$, $\lambda\in\Bbb T$.  
By (4.5) we have 
$$
V_\lambda B V_\lambda^*=\lambda B, 
$$
and it follows that both $\ker B$ and $\ker B^*$ are invariant 
under the action of $V$.  Since the spectral subspaces of $U$ and $V$
$$
H_k=\{\xi\in H: U_\lambda\xi=\lambda^k\xi,\lambda\in\Bbb T\},\quad
\tilde H_k=\{\zeta\in\tilde H: V_\lambda\zeta=\lambda^k\zeta,\lambda\in\Bbb T\}
$$
are related by $\tilde H_k = H_k\otimes\Lambda\Bbb C^d$, it follows 
that both $\ker B$ and $\ker B^*$ decompose into orthogonal sums 
$$
\ker B=\sum_k \ker B\cap \tilde H_k, \quad 
\ker B^*=\sum_k \ker B^*\cap \tilde H_k.  
$$
We conclude that 
$$
\ker D = \sum_k\ker D\cap\tilde H_k=\sum_k\ker B\cap\ker B^*\cap\tilde H_k.  
$$

It was shown in \cite{3}, Proposition 5.4, that each $H_k$ is a 
finite dimensional subspace of $M_H$, hence $\tilde H_k$ is a 
finite dimensional subspace of $M_H\otimes\Lambda\Bbb C^d$.  Since the 
restriction $B_{M_H}$ of $B$ to $M_H\otimes\Lambda\Bbb C^d$ 
is the boundary operator 
of the Koszul complex of $M_H$ it follows that 
for the restriction $B_k$ of $B$ to $M_H\cap\tilde H_k$ 
we have 
$$
\dim(\ker D\cap\tilde H_k)=\dim(\ker B\cap\ker B^*\cap\tilde H_k)=
\dim(\ker B_k/{\text{ran }}B_{k-1}).  
$$
Summing on $k$ we find that 
$$
\dim\ker D = \dim(\ker B_{M_H}/{\text{ran\,}}B_{M_H}).  
$$
The right side of the preceding formula is finite, because 
the Koszul complex of $M_H$ has finite dimensional
cohomology by Lemma 2.  

By restricting this argument respectively to the even and 
odd subspaces of $\tilde H$, one finds in the same way that
$\dim\ker D_+$ and $\dim\ker D_+^*$ are, respectively, the total 
dimensions of the even and odd cohomology of the Koszul 
complex of $M_H$, 
and that gives the two formulas of (4.6).  
\qedd\enddemo

\remark{Concluding remarks and examples}
It is natural to ask if Theorem B remains valid when one 
drops the hypothesis that $\bar T$ is graded.  On the surface, 
this may appear a foolish question since it is not known if 
the Dirac operator associated with a finite rank 
pure $d$-contraction is Fredholm, and if it is not Fredholm then
what does 
the index of $D_+$ mean?   The Dirac operator is known to be 
Fredholm for a class of concrete examples (this is an unpublished 
result of the author's), but the 
issue of Fredholmness for general pure finite rank $d$-contractions 
remains somewhat mysterious.  

Nevertheless, Stephen Parrott has proved a result \cite{11} for 
single operators that implies that Theorem B 
is true {\it verbatim} for the one-dimensional 
case $d=1$ and an arbitrary pure 
contraction $T$ of finite rank, graded or not.  His result implies 
that $T$ is necessarily a Fredholm operator.  
\endremark

Finally, we describe a class of examples of 
finite rank pure $d$-contractions $\bar T=(T_1,\dots,T_d)$ in 
arbitrary dimension $d=1,2,\dots$.  Most are ungraded.  
We show that these examples are Fredholm and we 
compute all three integer invariants 
(the index of the Dirac operator, the curvature invariant $K(\bar T)$, and 
the Euler characteristic $\chi(\bar T)$ of \cite{3}).  For some of these 
examples the formula $K(\bar T)= \chi(\bar T)$ of 
(\cite{3}, Theorem B) holds, but for most of them it fails.  
On the hand, in all cases the formula of Theorem B above 
$$
(-1)^dK(\bar T)=\dim\ker D_+ -\dim\ker D_+^*  \tag{4.7}
$$
is satisfied.  Indeed, we know of no examples for which 
(4.7) fails.  

Fix $d=1,2,\dots$ and let $r$ be a positive integer.  Following 
the notation and terminology of \cite{2}, \cite{3} we will 
consider the $d$-shift $\bar S=(S_1,\dots,S_d)$ of multiplicity 
$r+1$.  $\bar S$ acts on the Hilbert space $(r+1)\cdot H^2$, a 
direct sum of $r+1$ copies of the basic free Hilbert module 
$H^2=H^2(\Bbb C^d)$.  We consider certain invariant subspaces 
$M\subseteq (r+1)\cdot H^2$ and their quotient Hilbert modules 
$H=(r+1)\cdot H^2/M$.  The $d$-shift compresses to a pure 
$d$-contraction $\bar T=(T_1,\dots,T_d)$ acting on $H$, and 
the rank of $\bar T$ is at most $r+1$.  For the examples below, 
the rank is $r+1$ and $\bar T$ will have the properties asserted
above.  
The subspaces $M$ are defined as follows.  
Let $\phi_1,\phi_2,\dots,\phi_r$ be 
a set of multipliers of $H^2$ and set 
$$
M=\{(f,\phi_1f,\phi_2f,\dots,\phi_rf): f\in H^2\}\subseteq 
(r+1)\cdot H^2.  
$$

\proclaim{Proposition} Assume that the set of 
$r+1$ functions $\{1,\phi_1,\phi_2,\dots,\phi_r\}$ is linearly 
independent, and let 
$\bar T=(T_1,\dots,T_d)$ be the $d$-tuple of operators 
associated with the quotient Hilbert module $H=(r+1)\cdot H^2/M$.  
$\bar T$ is a pure $d$-contraction of rank $r+1$, it is Fredholm, 
and its index and curvature invariant are given by 
$$
\dim\ker(D_+)-\dim\ker D_+^*=(-1)^d\cdot r, \qquad K(\bar T) = r.  
$$
If each $\phi_k$ is a homogeneous polynomial of some degree
$n_k$ then the Euler characteristic is also given by 
$\chi(\bar T) = r$.  On the other hand, if 
$M$ contains no nonzero element $(p_1,p_2,\dots,p_{r+1})$ 
with polynomial components $p_k$, then $\chi(\bar T)=r+1$.  
\endproclaim

\remark{Remark}
For example, if each $\phi_k$ is the exponential of some nontrivial 
polynomial (or more generally, if no $\phi_k$ is a rational 
function), then the only $r+1$-tuple of polynomials 
$(p_0,p_1,\dots,p_r)$  that belongs to $M$ is the zero $r+1$-tuple.  
\endremark

\demo{proof}
We merely sketch the key elements of the proof.  

We deal first with the Euler characteristic.  This invariant is 
associated with the finitely generalted 
algebraic $\Bbb C[z_1,\dots,z_d]$-module 
$$
M_H={\text{span}}\{f(T_1,\dots,T_d)\zeta: f\in\Bbb C[z_1,\dots,z_d], 
\quad\zeta\in\Delta H\},
$$
$\Delta$ being the finite rank defect operator 
$\Delta=(\bold 1-T_1T_1^*-\dots-T_dT_d^*)^{1/2}$.  
Realizing the 
quotient $H=(r+1)\cdot H^2/M$ as the orthogonal complement 
$M^\perp\subseteq (r+1)\cdot H^2$, let $E_0\in\Cal B((r+1)\cdot H^2)$
be the projection onto the $r+1$-dimensional 
space of constant vector functions.  
The operators $T_1,\dots,T_d$ are obtained by compressing $S_1,\dots,S_d$ 
to $M^\perp$, and a straghtforward computation shows that 
$\Delta$ is identified with 
the square root of the compression of $E_0$ to $M^\perp$.  This operator 
is of rank $r+1$ because of the linear independence hypothesis on 
$\{1,\phi_1,\dots,\phi_r\}$ (for example, see 8.4.3 of \cite{2}).  
It follows that $M_H$ is identified with the projection onto $M^\perp$ 
of the space of all vector polynomials 
$$
S = \{(p_0,p_1,\dots,p_r): p_k\in\Bbb C[z_1,\dots,z_d]\}.  
$$
The $\Bbb C[z_1,\dots,z_d]$-module action of a polynomial $f\in\Bbb C[z_1,\dots,z_d]$
on $M_H$ is given by 
$$
f\cdot P_M^\perp(p_0,p_1,\dots,p_r) = P_M^\perp(fp_0, fp_1,\dots,fp_r).  
$$
Now assume that $M$ contains no nonzero element having polynomial 
components, and let $F= (r+1)\cdot\Bbb C[z_1,\dots,z_d]$ denote 
the free $\Bbb C[z_1,\dots,z_d]$-module of rank $r+1$.  
Consider the linear map 
$$
L:(f_0,f_1,\dots,f_r)\in F\mapsto 
P_{M^\perp}(f_0,f_1,\dots,f_r)\in M_H.  
$$
$L$ is injective by hypothesis, its range is all of $M_H$, 
and it is obviously a homomorphism of $\Bbb C[z_1,\dots,z_d]$-modules.  
Hence $M_H$ is a free module of rank $r+1$ and its Euler 
characteristic is $r+1$.  This shows that $\chi(\bar T)=r+1$ in this 
case.  

On the other hand, if each $\phi_k$ is a homogeneous polynomial then 
one may extend the argument of the proof of 
(\cite{2}, Proposition 7.4) (which addresses the 
case $r=1$ explicitly) in straightforward way to 
show that $\bar T$ is a graded $d$-contraction.  It follows from 
 Theorem B of \cite{3} that $\chi(\bar T) = K(\bar T)$.  We will 
show momentarily that in all cases we have $K(\bar T)=r$, and this 
calculates the Euler characteristic for the asserted cases.  

We show next  
that $\bar T$ is Fredholm of index $(-1)^dr$.  For that, 
it is enough to show that $\bar T$ is similar 
to a Fredholm $d$-tuple whose index is known to be $(-1)^dr$.  
The latter $d$-tuple is the $d$-shift of multiplicity $r$.  In more
detail, consider the linear mapping $A:(r+1)\cdot H^2\to r\cdot H^2$ 
defined by 
$$
A(f_0,f_1,f_2,\dots,f_r) = (f_1-\phi_1f_0,f_2-\phi_2f_0,\dots,f_r-\phi_rf_0).  
$$
It is clear that $A$ is bounded, surjective, has kernel $M$, and 
intertwines the action of $\bar S$ (acting on $(r+1)\cdot H^2$) and 
the multiplicity $r$ $d$-shift acting on $r\cdot H^2$.  Thus 
$A$ promotes to an isomorphism of Hilbert spaces 
$\tilde A: H\to r\cdot H^2$ which implements a similarity of 
$\bar T$ and the $d$-shift of multiplicity $r$.  The latter 
is known to be Fredholm and has index $(-1)^dr$.  

Finally, in order to calculate 
$K(\bar T)$ we appeal to a result of 
Greene, Richter and Sundberg \cite{9} as follows.  
Identifying $H$ with $M^\perp\subseteq (r+1)\cdot H^2$, 
we have already seen that 
the natural projection $L=P_{M^\perp}: (r+1)\cdot H^2\to H$ is the minimal 
dilation of $H$ in the sense of \cite{2}, and obviously 
$L$ is a co-isometry with $L^*L=\bold 1-P_M$.  Now if one 
evaluates all of the functions in $M$ at a point $z$ in the open unit ball 
of $\Bbb C^d$, one obtains the following linear subspace of 
$\Bbb C^{r+1}$
$$
M(z) = \{(\lambda,\lambda\phi_1(z),\lambda\phi_2(z),\dots,\lambda\phi_r(z)): 
\lambda\in\Bbb C\}.  
$$
This is 
a one-dimensional space having codimension $(r+1)-1=r$, and the 
same assertion is valid for almost every point $z$ on the 
boundary of the unit ball.  By the results of \cite{9}, the
codimension of $M(z)\subseteq \Bbb C^{r+1}$ is equal to 
$K(\bar T)$ for almost every $z$ 
in the boundary of the unit ball.  Thus, $K(\bar T)=r$.  
\qedd\enddemo

\end